%% file: main_amsart.tex
\documentclass{amsart}

\usepackage[pagewise]{lineno}

\usepackage{graphicx}
\usepackage{amsmath}
\usepackage{amsthm}
\usepackage{amsfonts}
\usepackage{amssymb}
\usepackage{mathtools}
   
\usepackage{amsaddr}

\usepackage{subfiles}

\usepackage[dvipsnames]{xcolor}

\usepackage[square,numbers]{natbib}

\usepackage{enumerate}
\usepackage[shortlabels]{enumitem}

\theoremstyle{plain}
\newtheorem{thm}{Theorem}[section]
\newtheorem{prop}{Proposition}[section]
\newtheorem{lem}{Lemma}[section]

\newtheorem{rem}{Remark}[section]
\newtheorem{que}{Question}[section]

\theoremstyle{definition}
\newtheorem{defi}{Definition}[section]

\allowdisplaybreaks

\title{Uniformly Expanding Full Branch Maps with a Wild Attracting Point}
\author[1]{Rubio Gunawan}
\address[1]{Scuola Internazionale Superiore di Studi Avanzati (SISSA), Trieste, Italy}
\address[1]{Abdus Salam International Centre for Theoretical Physics (ICTP), Trieste, Italy.}
\thanks{The author would like to thank Stefano Luzzatto for his supervision, and Hamza Ounesli for the discussion on the result and consequences.}

\date{\today}

\begin{document}

\begin{abstract}
    We give a counterintuitive example of a uniformly expanding full branch map 
    such that the point zero is a wild attractor.
\end{abstract}

\maketitle

\section{Introduction and Main Result}
\subfile{Sections/Intro}

\section{Outline of Proof}
\subfile{Sections/Proof}

\section{Topology of Basin}
\subfile{Sections/Topology}

\section{Invariant Subsets and Cylinder Sets}
\subfile{Sections/Subset}

\section{Convex Preimages of $C$}
\subfile{Sections/Convex}

\section{Full Measure Basin}
\subfile{Sections/Full}

\section{Construction of an Example}
\subfile{Sections/Existence}

\bibliographystyle{unsrt}
\bibliography{ref}

\end{document}

%% file: Sections/Intro.tex
\subsection{Definition and Statement of Main Theorem}

Let \( I \) denote a compact interval, \( \mathcal P \) a partition of \( I \) (mod 0) into subintervals, and \( r \geq 1\).  
We recall first of all the following  standard definition.

\begin{defi}
 \( f: I \to I \) is a  \emph{piecewise \( C^{r}\), uniformly expanding,  full branch map}  with respect to \( \mathcal P \) if there exists a constant \( \lambda > 1 \) such that for every \( J \in \mathcal P\) the map \( f: \mathring J \to \mathring I \) is a \( C^{r} \) diffeomorphism and \( |f'(x)|\geq \lambda \) for all \( x\in \mathring J \). 
\end{defi}

By classical results, if  \( r > 1 \)  and  \( \mathcal P \) is finite then \( f \) admits an \emph{ergodic invariant probability measure equivalent to Lebesgue}. When \( \mathcal P \) is infinite, the same result holds but under some additional, albeit fairly natural, conditions such as Adler's condition \cite{Adl73} on the ratio \( |f''(x)|/|f'(y)|^2 \) between the first and second derivatives
or Rychlik's condition  \cite{Ryc76} on  the variation of \( 1/|f'|\). Not much attention has been paid however  to the question of  what can happen if such conditions fail.  In this paper we will prove the following remarkable and counter-intuitive fact.

\begin{thm}
\label{Thm_Main}
    There exists (many) $f: I \to I$  piecewise analytic uniformly expanding full branch map with a wild fixed point attractor with full measure basin. 
\end{thm}

By the basin of a point \( p \) we mean the set 
\begin{equation}\label{eq:basinp}
B_p:= \{x| f^n(x) \xrightarrow[]{n \to \infty} p\}.
\end{equation}
of points whose orbits actually converge to \( p \). Following Milnor \cite{Mil99}, we say that \( p \) is a \emph{wild} attractor if its basin has positive Lebesgue measure but is topologically \emph{meagre}, in particular its basin does not contain a neighbourhood of the fixed point. For this reason we cannot say that \( p \) is an ``attracting fixed point'', since this is usually interpreted to mean that its basin contains a neighbourhood of \( p \). The mechanism which causes \( p \) to be attracting in our examples is very different from the usual mechanism.
Instead of having an invariant neighborhood that contracts to $p$, we have a collection of disjoint intervals that each map to a neighborhood of $p$. Said intervals accumulate at $p$, with their Lebesgue density approaching $1$.
In fact part of what makes this result so counterintuitive is that the map is actually strongly ``expanding'' and also has positive (in fact infinite) topological entropy. It is nevertheless, as will be clear from the construction, quite a robust mechanism. 

\subsection{Background and Literature}
There seem to be just a handful of papers in the literature exploring these kinds of ``pathological'' examples, even though they are clearly very interesting. 

One approach is to relax the regularity condition, i.e. to relax the  \( r >1 \) condition and  consider \( C^{1} \) uniformly expanding full branch maps. It turns out that this can produce new phenomena which do not exist for \( r >1 \). 

A relatively explicit example of this kind was constructed in \cite{GS89} with finitely many branches and without any invariant probability measure absolutely continuous with respect to Lebesgue (\emph{acip}). 
They used an induced map and showed that, whether or not the induced map has an acip, their example cannot have an acip. However, they do not actually show if the induced map has an acip or not, so not much is known about the explicit dynamics of their example.
\begin{que}
    Does the map studied in \cite{GS89} have a (singular) physical measure? What is the support, and what are the dynamical properties of the measure?
\end{que}
In \cite{Gun26}, we apply the results of this paper to find singular physical measures.
We conjecture that a modified version of the arguments within our two papers may be able to describe the dynamics of \cite{GS89} more explicitly.

It was proved in 
\cite{CQ01} that in fact \emph{generic} $C^1$ expanding circle maps (which can be seen as special cases of piecewise \( C^{1}\) uniformly expanding full branch interval maps with finitely many branches) \emph{do not} admit any \emph{acip}. 
Instead they show that a generic map has a unique physical measure which is singular and fully supported. 
More explicit properties of these measures are still not known and would be well worth investigating. 
\begin{que}
    Which singular measures can be realized as a physical measure of an expanding circle map?
    Given any compact set $\Lambda$, is there an expanding circle map with a physical measure on $\Lambda$? Is it possible to produce an explicit example?
\end{que}

Closer in spirit to our result are the examples in \cite{Bug85} of uniformly expanding maps with infinite branches that do not admit an \emph{acip}. But these examples are crucially \emph{not full branch}, in fact the infimum of the length of the images is $0$, and so the mechanism is very different from what we present here. 
A similar kind of map was also considered in \cite{KN95} as a simplified model of the induced maps of Fibonacci maps.
Particular fractal and thermodynamic properties of this particular model are further studied in \cite{SV97, BT12}.

The existence of wild \emph{Cantor} attractors (rather than fixed point attractors) was  proved and studied in \cite{BKNS96, MS14} for a class of unimodal maps with a very special Fibonacci combinatorial structure. These maps are however not uniformly expanding and the proof uses completely different arguments and techniques based on the combinatorics of the orbit of the critical point. 

%% file: Sections/Proof.tex
\label{Section_Result}
For simplicity we let \( I = [0,1]\) and 
let  $\mathcal{F}$ denote  the set of full branch maps $f: [0,1] \to [0,1]$ with a countable generating partition
\(
\mathcal{P} = \{I_{n}\}_{n\in \mathbb N}
\)
which  accumulates (monotonically in \( n \)) at $0$ and such that each branch \( f|_{I_{n}}: I_{n}\to (0,1) \) is a non-singular orientation-preserving homeomorphism.  
For any \( f\in \mathcal F \) we define the following notation. 
For  \( n\in \mathbb N \) we split the interval \( [0,1]\) into a left and right part
 \[
 I_{n}^{-}\coloneq \bigcup_{k \geq n+1} I_k 
 \quad \text{ and } \quad 
 I_{n}^{+}\coloneq\bigcup_{k \leq n} I_k 
 \]
where \(  I_{n}^{-}\ \) contains all \( I_{k}\) to the left of \( I_{n}\) and \(  I_{n}^{+}\ \)  all those to the right of \( I_{n}\) as well as  \( I_{n}\) itself. Since each \( f|_{I_{n}}: I_{n}\to (0,1) \)  is an orientation-preserving  homeomorphism, we can write \( I_{n} = L_{n}\cup R_{n}\) as the union of two subintervals 
\[
L_{n}\coloneq \{x\in I_{n}: f(x) \in  I_{n}^{-}\}
\quad\text{ and } \quad
R_{n}\coloneq \{x\in I_{n}: f(x) \in  I_{n}^{+}\}.
\]
Using this notation we define three subfamilies of $\mathcal{F}$:
\begin{itemize}
    \item $\mathcal{F}_* \subset \mathcal{F}$ is the subfamily of maps such that there exists a monotonically increasing sequence $\{p_n\}_{n \geq 1}$ of positive numbers such that:
        \begin{equation}
        \label{eqn_F_star}
             \frac{|L_n|}{|I_n|} > p_{n} 
             \quad \text{ and } \quad 
             P:=\prod_{n \geq 1} p_n > 0.  
        \end{equation}
    This condition essentially says that $|L_n|/|I_n| \to 1$ fast enough, and a large subset of points (at least of measure $P$) keeps moving to the left.
    \item $\mathcal{F}_*^{weak} \subset \mathcal{F}_*$ is the subfamily of maps such that each restriction $f|_{L_n}:L_n \to I_{n}^{-}$ is convex. 
    \item $\mathcal{F}_*^{strong} \subset \mathcal{F}_*^{weak}$ is the subfamily of maps such that each branch $f|_{I_n}:I_n \to I$ is convex.
\end{itemize}

One should expect that for $f \in \mathcal{F}_*$, 
most orbits would tend to move to the left.
In fact when we consider a simplified case where $f \in \mathcal{F}_*$ is piecewise affine, the following remark shows we have transient Markov Chain behavior.
\begin{rem}
    \label{rem_transient}
    Suppose $f \in \mathcal{F}_*$ and has affine pieces $f|_{L_n}$ and $f|_{R_n}$ for all $n$. $f$ is equivalent to a Markov Chain with state space $\{1,2,3,...\}$.
    Because $P > 0$, the Markov Chain is transient.
\end{rem}

Letting \( B_{0}\) denote the basin of attraction of 0, recall \eqref{eq:basinp}, 
the technical core of the paper is to use the convexity of $\mathcal F_{*}^{weak}$ and $\mathcal F_{*}^{strong}$ 
to extend the observation of Remark \ref{rem_transient} 
to then obtain lower bounds for $|B_0|$.

\begin{prop}
\label{Prop_Positive}
For every $f \in \mathcal F_{*}^{weak}$, we have \(|B_0| > 0\).
\end{prop}

\begin{prop}
\label{Prop_Full}
For every $f \in \mathcal F_{*}^{strong}$, we have \(|B_0| = 1\).
\end{prop}

\begin{figure}[ht]
\centering
\caption{Example of a map in \(\mathcal F_{*}^{strong}\).}
\includegraphics[width=0.4\textwidth]{"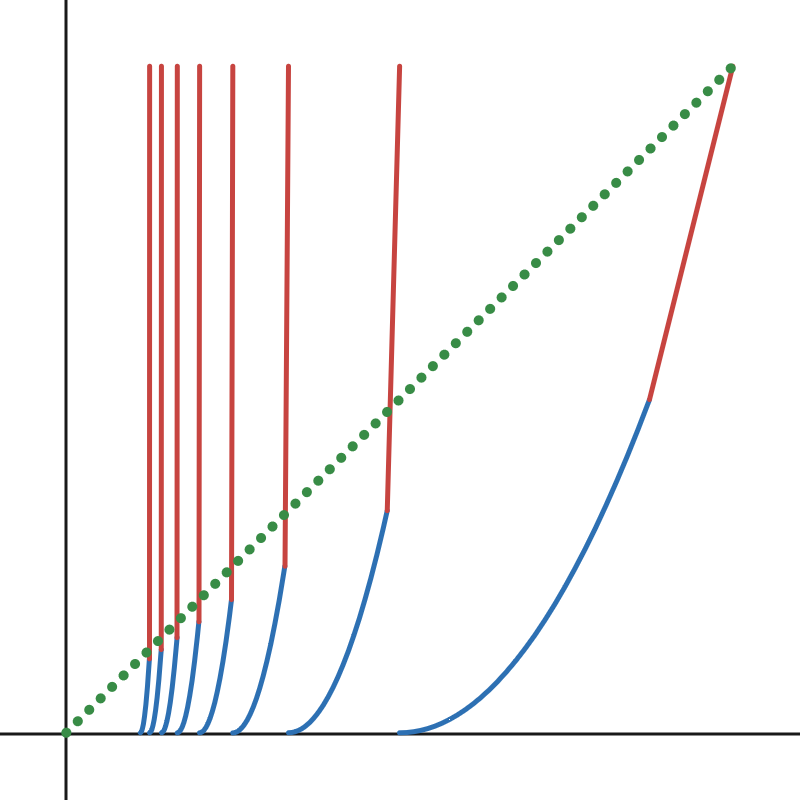"}
\end{figure}

Notice that \(\mathcal F_{*}^{strong}\) only requires the branches $f|_{I_n}$ to be convex non-singular orientation-preserving \emph{homeomorphisms}. It is not immediately obvious that the conditions which define \( \mathcal F_{*}^{strong}\) are compatible with each branch being analytic and uniformly expanding. Thus we prove the following existence result.
\begin{prop}[]
\label{Prop_Existence}
    There exists (many) $f \in \mathcal F_{*}^{strong}$, uniformly expanding and piecewise analytic. 
\end{prop}
Finally, we also have this relatively simple topological result.
\begin{prop}
\label{Prop_Meagre}
    Let $f \in \mathcal F$.
    Then $B_0$ is a meagre set.
\end{prop}
It is easy to see that the Propositions \ref{Prop_Full}, \ref{Prop_Existence}, and \ref{Prop_Meagre} together imply Theorem \ref{Thm_Main}.
We briefly outline the contents of the remaining sections.
In Section \ref{Section_Topology} we will prove Proposition \ref{Prop_Meagre}.
Section \ref{Section_Subset} presents particular invariant subsets of $[0,1]$ that play a role in the dynamics of $f \in \mathcal{F}$. In Section \ref{Section_Convex} we prove the technical Proposition \ref{Prop_Tech} that gives a lower bound for the Lebesgue measure of those subsets, to prove Proposition \ref{Prop_Positive}. In Section \ref{Section_Full} we apply that technical proposition for $f \in \mathcal{F}_*^{strong}$, to prove Proposition \ref{Prop_Full}.
Section \ref{Section_Existence} contains an explicit construction of a map $f \in \mathcal{F}_*^{strong} $ that proves Proposition \ref{Prop_Existence}.

%% file: Sections/Topology.tex
\label{Section_Topology}
In this section we prove Proposition \ref{Prop_Meagre}.
Let $f\in \mathcal F$.
Then $f$ has a generating partition, thus it is locally eventually onto.
We now show that $B_0$ is meagre.
First we define the following sets:
\begin{align*}
  M_n := \{x|f^k(x) \notin I_1, \forall k \geq n\},\quad
  M := \bigcup_{n\geq 1} M_n.
\end{align*}
$M_n$ is the set of points that do not recur to $I_1$ after time $n$,
and $M$ is the set of points that eventually do not recur to $I_1$.
Because $I_1$ is not a neighborhood of $0$,
points that converge to $0$ must eventually not recur to $I_1$.
Therefore, $B_0 \subset M$.

\begin{lem}
    \label{Lem_Nonrecur}    
    The set $M$ is meagre.
\end{lem}
\begin{proof}
    Fix any index $n \geq 1$. Fix any open interval $J$.
    Because $f$ is locally eventually onto, there must be a sufficiently large $m \geq n$
    and a subinterval $J' \subset J$ that is homeomorphic by $f^m$ to $I_1$.
    Thus $\overline{M_n}$ is not dense in $J$, for any open interval $J$. 
    Therefore $\overline{M_n}$ is nowhere dense, for any index $n \geq 1$.
    Thus the countable union $M = \bigcup_{n\geq 1} M_n$ is a meagre set.
\end{proof}

\begin{proof}[Proof of Proposition \ref{Prop_Meagre}]
    Because $B_0 \subset M$, by Lemma \ref{Lem_Nonrecur}, $B_0$ is also meagre.
\end{proof}

%% file: Sections/Subset.tex
\label{Section_Subset}
Let $f \in \mathcal F$.
We introduce two particular subsets of $I$:

\begin{align*}
C &:= \left\{x \in I: f^i(x) \in \bigcup_{j \geq 0} L_j \text{ for all } i \geq 0\right\}.\\
E &:= \left\{x \in I: \exists N \geq 0, f^i(x) \in \bigcup_{j \geq 0} L_j \text{ for all } i \geq N\right\}.
\end{align*}
Then $C$ is forward invariant, and it contains points that have an orbit that stays inside the left subintervals $L_n$, which means they keep moving to the left. 
Because the intervals $I_n$ accumulate at $0$ on the left, $C \subset B_0$.
Moreover $E$ is fully invariant, and it contains points that eventually fall into $C$. 
Thus, $C \subset E \subset B_0$.
We are naturally interested in studying the structure and measure of $C$ and $E$.

To begin the study of $C$, we define a family of sets indexed by positive integers. 
Given $n \geq 1$, define
\[
C_n := \left\{x \in I: f^i(x) \in \bigcup_{j \geq 0} L_j \text{ for all } 0 \leq i < n \right\}.
\]
Immediately by the definition, $C = \bigcap_{n \geq 1} C_n$.
To study the structure of $C_n$, note that, for $n=1$:
\[C_1 =  \left\{x \in I: x \in \bigcup_{j \geq 0} L_j\right\} = \bigcup_{j \geq 0} L_j.\]
For a more general description, we need to study the partitions of full branch maps and its refinements.
First, recall that $f$ is associated with a partition (modulo Lebesgue measure) of $I$ into subintervals given by $\mathcal{P} = \{I_1,I_2,I_3,...\}$.

Recall $f^{n+1}$ is also a full branch map for each $n \geq 1$, 
associated to a refined partition $\mathcal{P}_{n+1} = \bigvee_{i=0}^{n}f^{-i} \mathcal P$.
For the rest of the paper, we borrow the terminology and notation of cylinder sets from symbolic dynamics to describe the partition elements of $\mathcal{P}_{n+1}$. Given $n$ positive integers $(k_i)_{i=0}^n$, we define its cylinder set:
\[
[k_0, k_1, ..., k_n] := (I_{k_0}) \cap f^{-1}(I_{k_1}) \cap ... \cap f^{-n}(I_{k_n})
\in \mathcal{P}_{n+1}.
\]
Then we obtain the elements of $\mathcal{P}_{n+1}$ as cylinder sets:
\[
\mathcal{P}_{n+1} = \{ [k_0,...,k_n]: k_0,...,k_n \ge 1 \}.
\]
Cylinder sets also give a partition of $C_n$.
\begin{lem} [Partition by Cylinder Sets]
    \label {Lem_Cylinder_Partition}
    For each $n \geq 1$,
    \[
    C_n = \bigcup_{1 \leq k_0 < k_1 < ... < k_n} [k_0, k_1, ..., k_n].
    \]
\end{lem}
\begin{proof}
    Take any cylinder set $[k_0, k_1, ..., k_n]$.
    Take any index $i \in \{0,1,...,n-1\}$.
    By definition, $f^i(x) \in \bigcup_{j \geq 0} L_j \iff k_i < k_{i+1}$.
    Thus $C_n$ is the union of cylinder sets that are associated to strictly increasing finite sequences $(k_i)_{0\leq i \leq n}$.
\end{proof}

\subsection{Left Subintervals of Cylinders}
Recall from Section \ref{Section_Result}, for a given interval $I_n \in \mathcal P$,
we defined its left subinterval $L_n \subset I_n$.
Now, given a cylinder set $[k_0,k_1, ... k_n]$, we define its left subinterval by the following union:
\[L([k_0,k_1,...,k_n]) := \bigcup_{k_{n+1} \geq k_n + 1}^\infty [k_0,k_1,...,k_n,k_{n+1}].\]
By the following Lemma \ref{Lem_Image}, this new definition is compatible with the old one.

\begin{lem}[Image of left subinterval]
\label{Lem_Image}
Take any cylinder set $[k_0,k_1,...,k_n] \in \mathcal{P}_{n+1}$. Then:
    \begin{align*}
        f^n: [k_0,k_1,...,k_n] &\to I_{k_n} 
        &\text{is a homeomorphism by definition, and}\\
        f^n: L([k_0,k_1,...,k_n]) &\to L_{k_n}
        &\text{is also a homeomorphism.}
    \end{align*}
\end{lem}
\begin{proof}
    Because a restriction of a homeomorphism is a homeomorphism, we only need to compute the image:
    \begin{align*}
        f^n(L([k_0,k_1,...,k_n]))
        &=\bigcup_{k_{n+1} \geq k_n + 1}^\infty f^n([k_0,k_1,...,k_n,k_{n+1}])\\
        &=\bigcup_{k_{n+1} \geq k_n + 1}^\infty [k_n,k_{n+1}]
        =L_{k_n}.
    \end{align*}
\end{proof}

These left subintervals are also useful to describe how $C_{n+1}$ is nested inside $C_n$.
\begin{lem}[Partition by Left Subintervals of Cylinder Sets]
\label{Lem_Partition}
    Take any $n \geq 1$. 
    Then we can partition $C_{n+1}$ by the left subintervals of cylinder sets of $C_n$:
    \begin{align*}
        C_{n+1} = \bigcup_{1\leq k_0<...<k_{n+1}} [k_0,k_1,...,k_n,k_{n+1}] = \bigcup_{1\leq k_0<...<k_n} L([k_0,k_1,...,k_n]).
    \end{align*}
\end{lem}
\begin{proof}
    The first equation is true by Lemma \ref{Lem_Cylinder_Partition}.
    The second equation is true by definition of $L([k_0,k_1,...,k_n])$.
\end{proof}

By Lemmas \ref{Lem_Image} and \ref{Lem_Partition}, we can now picture how $C_{n+1}$ is nested in $C_n$.
\begin{figure}[ht]
    \centering
    \caption{Nested sets $I \supset C_1 \supset C_2$.}
    \includegraphics[width=0.7\linewidth]{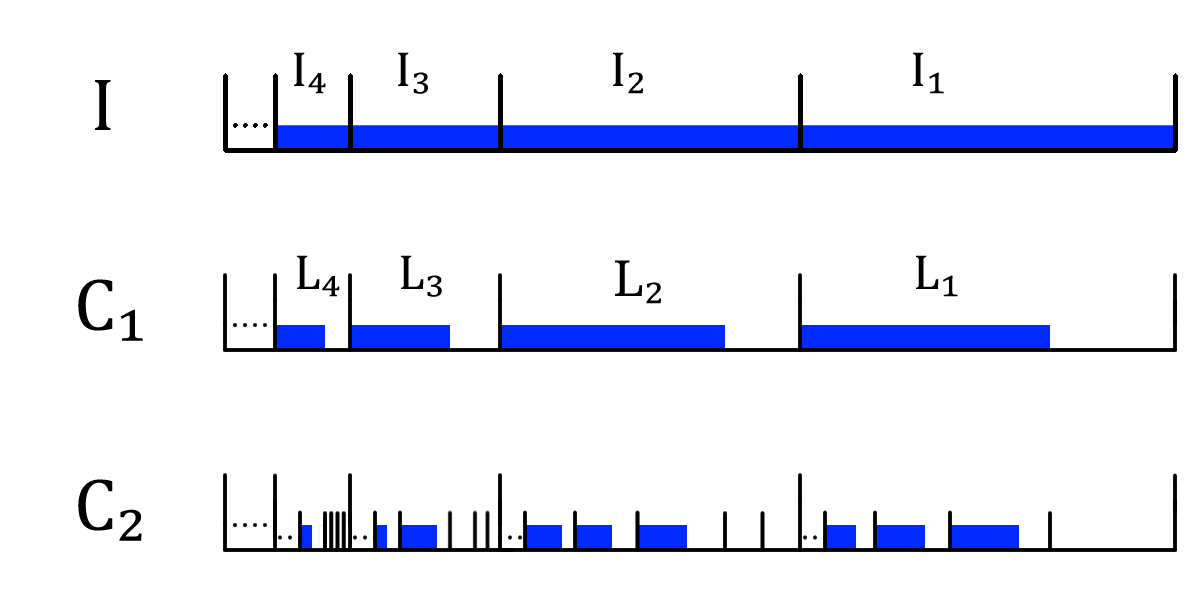}
\end{figure}

%% file: Sections/Convex.tex
\label{Section_Convex}
In the previous Section \ref{Section_Subset}, given a full branch map $f\in \mathcal F$, we have defined the nested sets $C \subset E \subset B_0$.
In this section we state and prove Proposition \ref{Prop_Tech},
which is the main technical tool for the results of this paper.

\begin{prop}(Convex Preimages of $C$)
\label{Prop_Tech}
    Let $f \in \mathcal F_*^{weak}$.
    Let $C$ be the set defined by $f$.
    Take any bounded interval $J$, and 
    $g : J \to (0,1]$ a convex orientation preserving homeomorphism. Then:
    \begin{align*}
        |g^{-1} (C_n)| \geq |J| \prod_{i=1}^n p_i, \forall n \geq 1.
        \quad \text{Furthermore, }|g^{-1} (C)| \geq |J| \prod_{i=1}^\infty p_i = |J| P.
    \end{align*}
\end{prop}
We begin with an elementary observation for preimages of left subintervals.
\begin{lem}[Convex Preimages of Left Subintervals]
\label{Lem_Convex}
    Let $J, I$ be bounded intervals.
    Let $g:J \to I$ be a convex orientation preserving homeomorphism.
    Let $L$ be a left subinterval of $I$ (i.e. they share a left endpoint.)
    Then the proportion of Lebesgue measure contained in the subinterval increases by each preimage:
\[\frac{|g^{-1}(L)|}{|J|} = \frac{|g^{-1}(L)|}{|g^{-1}(I)|} \geq \frac{|L|}{|I|}.\]
\end{lem}
\begin{figure}[ht]    
    \centering    
    \caption{Visual Proof of Lemma \ref{Lem_Convex}}    
    \includegraphics[width=0.4\textwidth]{"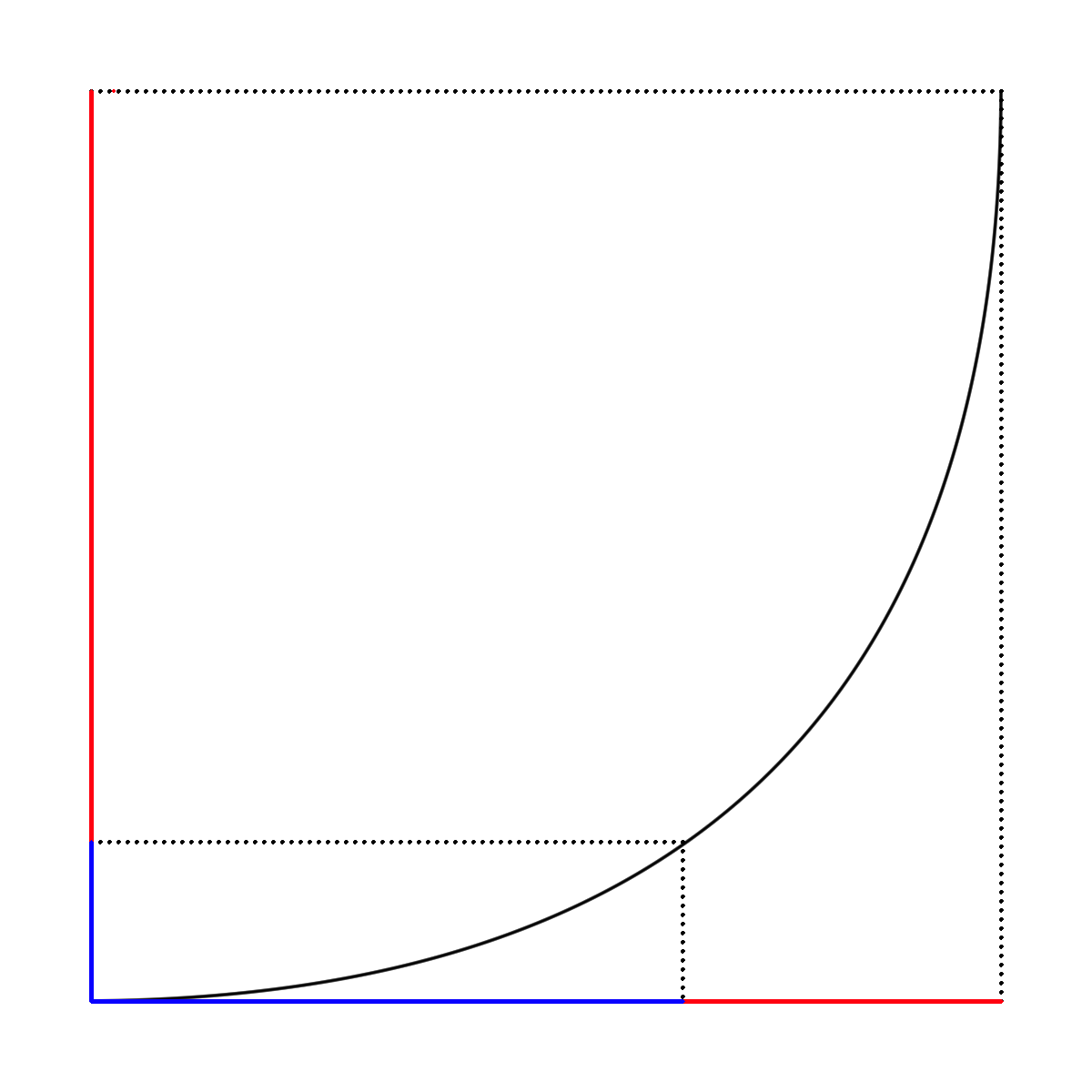"}
\end{figure}    

\begin{proof}[Proof of Proposition \ref{Prop_Tech}.]
    Our strategy is to apply Lemma \ref{Lem_Convex} with a proof by induction.
    For the base case $C_1$, recall that $C_1 = \bigcup_{n\geq 1}L_n$.
    We use this decomposition to obtain the estimate:
    \begin{align*}
        &|g^{-1}(C_1)|
        = \sum_{n \geq 1} |g^{-1}(L_n)|\\
        &\geq \sum_{n \geq 1} |g^{-1}(I_n)| \frac{|L_n|}{|I_n|}
        &\text{($g$ is convex, apply Lemma \ref{Lem_Convex})}\\
        &\geq \sum_{n \geq 1} |g^{-1}(I_n)| p_n
        &\text{(because $f \in \mathcal F$)}\\
        &\geq \left(\sum_{n \geq 1} |g^{-1}(I_n)|\right) p_1 = |J| p_1.
        &\text{($\{p_n\}_{n\geq1}$ is increasing)}\\
    \end{align*}
    Now we perform the inductive step.
    Recall the partition for $C_n$ by Lemma \ref{Lem_Cylinder_Partition}, 
    and the partition for $C_{n+1}$ by Lemma \ref{Lem_Partition}:
    \[
    C_n = \bigcup_{k_0<...<k_n} [k_0,...,k_n], \quad C_{n+1} = \bigcup_{k_0<...<k_n} L([k_0,...,k_n]).
    \]
    We use the two partitions to obtain the estimate:
    \begin{align*}
        &|g^{-1}(C_{n+1})|
        = \sum_{k_0<...<k_n} |g^{-1}(L([k_0,...,k_n]))|\\
        &\geq \sum_{k_0<...<k_n} |g^{-1}([k_0,...,k_n])| \frac{|L([k_0,...,k_n])|}{|[k_0,...,k_n]|}
        &\text{($g$ is convex, apply Lemma \ref{Lem_Convex})}\\
        &\geq \sum_{k_0<...<k_n} |g^{-1}([k_0,...,k_n])| \frac{|f^n(L([k_0,...,k_n]))|}{|f^n([k_0,...,k_n])|}
        &\text{($f^n$ is convex on $[k_0,...,k_n]$)}\\
        &= \sum_{k_0<...<k_n} |g^{-1}([k_0,...,k_n])| \frac{|L_{k_n}|}{|I_{k_n}|}
        &\text{(by Lemma \ref{Lem_Image})}\\
        &\geq \sum_{k_0<...<k_n} |g^{-1}([k_0,...,k_n])| p_{k_n}
        &\text{(because $f \in \mathcal F$)}\\
        &\geq \left(\sum_{k_0<...<k_n} |g^{-1}([k_0,...,k_n])|\right) p_{n+1}
        &\text{($\{p_n\}_{n\geq1}$ is increasing)}\\
        &= |g^{-1}(C_n)| p_{n+1} \geq |J| \prod_{i=1}^{n+1} p_i.
    \end{align*}
\end{proof}

\begin{proof}[Proof of Proposition \ref{Prop_Positive}]
Let $f \in \mathcal{F}_*^{weak}$. 
We apply Proposition \ref{Prop_Tech}, let $J = \enspace (0,1]$ and $g$ be the identity. 
By (\ref{eqn_F_star}), we obtain $|B_0| \geq |C| \geq |(0,1]| P > 0$.
\end{proof}

%% file: Sections/Full.tex
\label{Section_Full}

We want to find sufficient conditions for $|B_0| = 1$.
Recall in the proof of Proposition \ref{Prop_Positive}, we use $|C|$ as a lower bound for $|B_0|$.
Note that $C_1$ is obviously not full measure, so  $C$ cannot have full measure.
To show $|B_0| = 1$, we need a way to measure the elements of $B_0$ outside of $C$.
Recall from Section \ref{Section_Subset} the set $E \subset B_0$, 
which is the set of points that eventually fall into $C$.
Now define $D$ to be the set of points that do not fall into $C$ (excluding $0$):
\[ D := \enspace (0,1] \setminus E.\]
We want to show $|E| = 1$, which is equivalent to showing $|D| = 0$.
Our strategy is to express $D$ as a nested intersection of decreasing sets $D_m$.
Define $D_1 := (0,1] \setminus C$.
Define $\mathcal{D}_1$ to be the collection of cylinder sets associated to finite sequences that strictly increase except at the last index.
\[\mathcal{D}_1 := \{[k_0,k_1,...,k_n]: k_0<k_1<k_2<...<k_{n-1},k_{n-1}\geq k_n\}.\]
Notice that this collection is pairwise disjoint, and it forms a partition for $D_1$:
\[D_1 = \bigcup_{J \in \mathcal{D}_1} J.\]

To obtain $D_2$, we want to remove the points of $D_1$ that fall into $C$.
Let $J \in \mathcal{D}_1$ be a cylinder set which is also an element of $\mathcal{P}_n$.
Because $f^{n}|_J : J \to (0,1]$ is a homeomorphism, we can partition $J$ into two disjoint subsets:
\[(0,1] = D_1 \bigsqcup C \implies J = (f^{n}|_J)^{-1}{(D_1)} \bigsqcup (f^{n}|_J)^{-1}{(C)}\]
We now define $D_2$ by the union of preimages of $D_1$:
\begin{align*}
    &D_2 := \bigcup_{J \in \mathcal{D}_1} (f^{n}|_J)^{-1}{(D_1)}\\.
    &\text{(Note that $n$ depends on $J$, but we omit this in the notation for brevity.)}
\end{align*}
Similar to $D_1$, $D_2$ is also a union of pairwise disjoint cylinder sets. Note that:
\begin{align*}
    D_2 &= \bigcup_{J \in \mathcal{D}_1} (f^{n}|_J)^{-1}(D_1)
    = \bigcup_{J \in \mathcal{D}_1} (f^{n}|_J)^{-1}(\bigcup_{J' \in \mathcal{D}_1} J')
    = \bigcup_{J,J' \in \mathcal{D}_1} (f^{n}|_J)^{-1}(J').
\end{align*}
We collect these new cylinder sets in the family $\mathcal{D}_2$, and obtain a partition for $D_2$.
\begin{align*}
    \mathcal{D}_2 := \{(f^{n}|_J)^{-1}(J') | J,J' \in \mathcal{D}_1\}. \quad
    D_2 = \bigcup_{J \in \mathcal{D}_2} J.
\end{align*}
We repeat this process inductively for any $m \geq 1$:
\begin{align*}
    \mathcal{D}_{m+1} 
    := \{(f^n|_J)^{-1}(J')|J \in \mathcal{D}_m, J' \in \mathcal{D}_1\}.\quad
    D_{m+1} := 
    \bigcup_{J \in \mathcal{D}_{m+1}} J.
\end{align*}

We now study an individual cylinder set from $\mathcal{D}_{m+1}$.
Let $J = [k_0,...,k_{n-1}] \in \mathcal{D}_m \cap \mathcal{P}_n$, 
and $J'=[k_0'...,k'_{n'-1}] \in \mathcal{D}_1 \cap \mathcal{P}_{n'}$.
We obtain the new cylinder set by their concatenation:
\[(f^{n}|_J)^{-1}(J') = [k_0,...,k_{n-1},k_0',...,k'_{n'-1}].\]
In general, a cylinder set of $\mathcal{D}_m$ 
is a finite concatenation of $m$ cylinder sets of $\mathcal{D}_1$.
Now take the infinite nested intersection $\bigcap_{m \geq 1} D_m$.
By definition of $D_m$, the points of this infinite nested intersection will never fall into $C$. Therefore,
\[D = \bigcap_{m \geq 1} D_m.\]

Now we estimate $|D|$.
Recall that $D_{m+1}$ is obtained from $D_m$ by removing convex preimages of $C$.
By Proposition \ref{Prop_Tech}, we obtain the following bound:

\begin{lem}
    \label{Lem_D_m}
    Let $f \in \mathcal F_*^{strong}$. Then
    $|D_{m}| \leq (1-P)^{m}, \forall m \geq 1$.
\end{lem}
\begin{proof}
    Observe that for $m = 1$,
    $|D_1| = 1 -|C| \leq 1-P$.
    Now we obtain the estimate for general $m$ by induction.
    Recall we obtain $D_{m+1}$ by refining the cylinder sets of $\mathcal{D}_m$.
    Let $J \in \mathcal{D}_m$.
    Let $n \geq 1$ be such that $J \in \mathcal{P}_n$.
    Since $f \in \mathcal{F}_*^{strong}$,
    $f^n:J \to ]0,1]$ is a convex orientation preserving homeomorphism,
    then by Proposition \ref{Prop_Tech} we have
    $|(f^n|_J)^{-1}(C)| \geq |J|P$.
    We perform the inductive step:    
    \begin{align*}
        &|D_{m+1}|
        = \sum_{J\in \mathcal{D}_{m}} \left| J\cap (f^n|_J)^{-1}(D_1) \right|
        = \sum_{J\in \mathcal{D}_{m}} (|J| - |(f^n|_J)^{-1}(C)|)\\
        &\leq \sum_{J\in \mathcal{D}_{m}} (|J| - |J|P)
        = \left(\sum_{J\in \mathcal{D}_{m}} |J|\right) (1-P)
        =|D_m| (1-P) \leq (1-P)^{m+1}.
    \end{align*}
\end{proof}

\begin{proof}[Proof of Proposition \ref{Prop_Full}]
As a consequence of Lemma \ref{Lem_D_m}:
$|D| = \lim_{m \to \infty} |D_m| \leq \lim_{m \to \infty} (1-P)^m = 0$.
Furthermore $|E| = 1 - |D| = 1$, and $|B_0| = |E| = 1$.
\end{proof}

%% file: Sections/Existence.tex
\label{Section_Existence}

In this section we will construct a map $f \in \mathcal{F}_*^{strong}$,
that satisfies Proposition \ref{Prop_Existence}.
We begin by fixing several parameters, the intervals $I_n$, and subintervals $L_n$.

\begin{itemize}
    \item Fix $c \in (1,2)$. 
    We then fix the sequence of endpoints $a_n := c^{-n+1}$,
    which then gives us the intervals $I_n := [a_{n+1},a_n] = [c^{-n},c^{-n+1}]$.
    \item We fix the constant of expansion $\lambda \:= \frac{1}{c-1} > 1$.
    \item We fix an increasing sequence $\{p_n\}_{n \geq 1}$ such that 
    $\prod_{n\geq 1} p_n > 0$. For example, $p_n = 1 - (n+1)^{-2}$.
    \item We assume $p_1$ is large enough such that $\frac{p_1}{1-p_1} > \lambda$.
\end{itemize}

Inside each interval $I_n$, we define the left subinterval $L_n$ such that it has length $|L_n| = p_n |I_n|$. This also defines the right subinterval $R_n$.
To make $f \in \mathcal{F}_*^{strong}$ and uniformly expanding,
it is sufficient to define each branch $f|_{I_n}$
such that they satisfy the following three properties:
\begin{enumerate}[start = 1, label = {(P\arabic*)}:]
    \item $f|_{I_n} : I_n \to (0,1]$ is a convex orientation-preserving homeomorphism.
    \item $f'|_{\mathring I_n} \geq \lambda$.
    \item The left subinterval $L_n$ must be mapped to 
    $I_n^- =\bigcup_{k \geq n+1} I_k = (0,a_{n+1}]$.
\end{enumerate}
We also want to make $f$ as regular as possible.
The following computational lemma will help us define $f|_{I_n}$ in a way that satisfies all three properties:
\begin{lem}
    \label{Lem_Exp_Uniform}
    For each $n \geq 1$, we have:
    $
    {|I_n^-|}/{|L_n|} = {\lambda}/{p_n}, \quad
    {|I_n^+|}/{|R_n|} > {\lambda}/{p_n}.
    $
\end{lem}
\begin{proof}
    We prove the first item.
    \begin{align*}
        \frac{|I_n^-|}{|L_n|}
        = \frac{|I_n^-|}{p_n |I_n|}
        = \frac{1}{p_n}\frac{a_{n+1}}{a_n - a_{n+1}}
        = \frac{1}{p_n}\frac{c^{-n}}{c\cdot c^{-n} - c^{-n}}
        = \frac{1}{p_n}\frac{1}{c - 1} = \frac{\lambda}{p_n}.
    \end{align*}
    Now we prove the second item.
    \begin{align*}
        \frac{|I_n^+|}{|R_n|}
        &= \frac{1-a_{n+1}}{(1-p_n) |I_n|}
        = \frac{1}{(1-p_n)}\frac{1-a_{n+1}}{a_n - a_{n+1}} 
        = \frac{1}{(1-p_n)}\frac{1-c^{-n}}{c\cdot c^{-n} - c^{-n}}\\
        &= \frac{1}{(1-p_n)}\frac{c^n - 1}{c - 1}
        \geq \frac{1}{(1-p_n)}\frac{c - 1}{c - 1} = \frac{1}{(1-p_n)}.
    \end{align*}
    Now recall by our choice of $\{p_n\}_{n \geq 1}$,
    \begin{align*}
        \frac{p_1}{1-p_1} > \lambda 
        \implies \frac{p_n}{1-p_n} > \lambda 
        \implies \frac{1}{1-p_n} > \frac{\lambda}{p_n}.
    \end{align*}
    Then we obtain ${|I_n^+|}/{|R_n|} \geq {1}/{(1-p_n)}> {\lambda}/{p_n}$.
\end{proof}

\begin{proof}[Proof of Proposition \ref{Prop_Existence}]
By Lemma \ref{Lem_Exp_Uniform}, we have $|I_n^+|/|R_n| > |I_n^-|/|L_n| > \lambda$.
We then define each $g|_{I_n}: I_n \to I$ to be a strongly convex 
($g|_{I_n}''(x) > m_n > 0$) orientation-preserving $C^\infty$ homeomorphism with
$g'(x) \geq \lambda$, such that $g(L_n) = I_n^-, g(R_n) = I_n^+$.
(Explicitly, one can define two pieces $G|_{L_n},G|_{R_n}$ by polynomials, then smoothly glue the two pieces by standard methods such as convolution.)
Thus we obtain $g \in \mathcal{F}_*^{strong}$ which is uniformly expanding
and piecewise $C^\infty$.

By Whitney Approximation Theorem, we approximate each $g_n$ by an analytic function $f_n$ with arbitrarily close derivatives. 
A sufficiently close zeroth derivative preserves (P3) (with a new sequence of subintervals $L_n$).
A sufficiently close first derivative preserves (P2) (with a new smaller $\lambda$).
Because $g_n$ is strongly convex, a sufficiently close second derivative preserves (P1).
By the three preserved properties we obtain $f \in \mathcal F_*^{strong}$ which is uniformly expanding and piecewise analytic.
\end{proof}

\begin{figure}[ht]
\centering
\caption{Sketch of uniformly expanding $f \in \mathcal{F_*}^{strong}$, with $c=4/3$.}
\includegraphics[width=0.5\textwidth]{"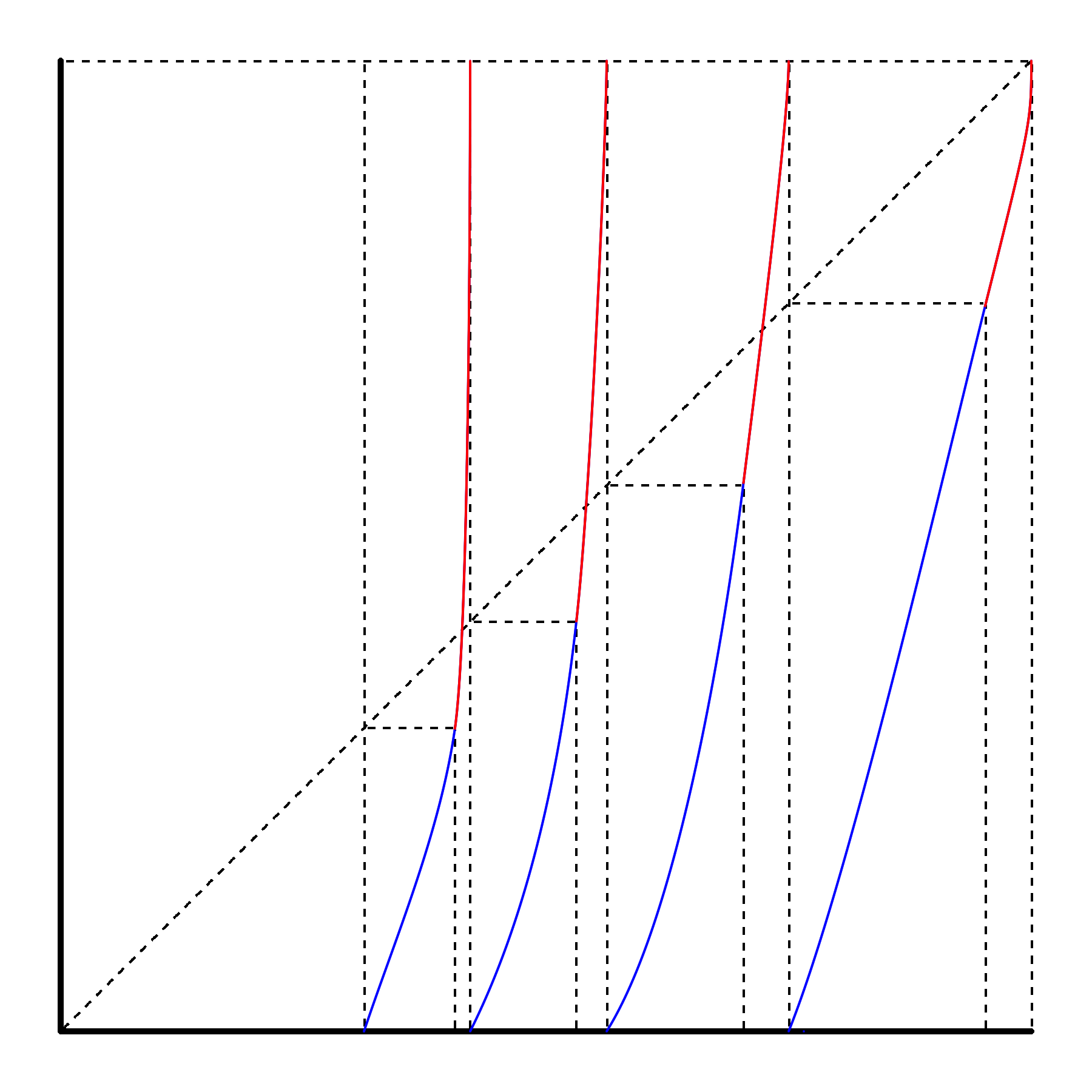"}
\end{figure}

\subsection{The Nonexpanding Case}
Recall we chose $c \in (1,2)$ in the beginning of this section.
By the following lemma, $c \in (1,2)$ is necessary to have uniform expansion.
\begin{lem}
    \label{Lem_Exp_Nonuniform}
    Fix $c \in [2,\infty)$.
    Let $f \in \mathcal{F}_*^{strong}$, 
    with branches given by $a_n = c^{-n}$.
    Then $f'(x)$ is not uniformly bounded away from $1$.
\end{lem}
\begin{proof}
    Because $f$ satisfies (A1), $\frac{|L_n|}{|I_n|}$ converges to $1$ as $n \to \infty$.
    We can now compute the limit of ${|f(L_n)|}/{|L_n|}$:
    \begin{align*}
    \lim_{n \to \infty}\frac{|f(L_n)|}{|L_n|}
    &= \lim_{n \to \infty}\frac{|I_n^-|}{|L_n|}
    = \lim_{n \to \infty}\frac{|a_{n+1}|}{|L_n|}\\
    &= \lim_{n \to \infty}\frac{a_{n+1}}{p_n|I_n|}
    = \lim_{n \to \infty}\frac{1}{p_n} \lim_{n \to \infty}\frac{a_{n+1}}{|I_n|}\\
    &= 1 \cdot \lim_{n \to \infty}\frac{a_{n+1}}{|I_n|}
    = \lim_{n \to \infty}\frac{a_{n+1}}{a_n - a_{n+1}}\\
    &= \lim_{n \to \infty}\frac{c^{-n}}{c^{-n+1} - c^{-n}}
    = \lim_{n \to \infty}\frac{c^{-n}}{c \cdot c^{-n} - c^{-n}}\\
    &= \frac{1}{c - 1}
    \leq \frac{1}{2 - 1} 
    = 1.
    \end{align*}
    Then for any $\varepsilon > 0$, there exists $n$ sufficiently large such that
    ${|f(L_n)|}/{|L_n|} \leq 1 + \varepsilon$.
    By the mean value theorem, there exists $x \in L_n$ such that $f'(x) \leq 1 + \varepsilon$,
    which means $f'(x)$ is not uniformly bounded away from $1$.
\end{proof}

%% file: ref.bib
@article {Mil99,
    AUTHOR = {Milnor, John},
     TITLE = {On the concept of attractor},
   JOURNAL = {Comm. Math. Phys.},
  FJOURNAL = {Communications in Mathematical Physics},
    VOLUME = {99},
      YEAR = {1985},
    NUMBER = {2},
     PAGES = {177--195},
      ISSN = {0010-3616,1432-0916},
   MRCLASS = {58F12 (58F08)},
  MRNUMBER = {790735},
MRREVIEWER = {Hans\ G.\ Bothe},
       URL = {http://projecteuclid.org/euclid.cmp/1103942677},
}

@article {BKNS96,
    AUTHOR = {Bruin, H. and Keller, G. and Nowicki, T. and van Strien, S.},
     TITLE = {Wild {C}antor attractors exist},
   JOURNAL = {Ann. of Math. (2)},
  FJOURNAL = {Annals of Mathematics. Second Series},
    VOLUME = {143},
      YEAR = {1996},
    NUMBER = {1},
     PAGES = {97--130},
      ISSN = {0003-486X,1939-8980},
   MRCLASS = {58F12 (30D05 58F23)},
  MRNUMBER = {1370759},
MRREVIEWER = {Hartje\ Kriete},
       DOI = {10.2307/2118654},
       URL = {https://doi.org/10.2307/2118654},
}

@article {GS89,
    AUTHOR = {G\'ora, P. and Schmitt, B.},
     TITLE = {Un exemple de transformation dilatante et {$C^1$} par morceaux
              de l'intervalle, sans probabilit\'e{} absolument continue
              invariante},
   JOURNAL = {Ergodic Theory Dynam. Systems},
  FJOURNAL = {Ergodic Theory and Dynamical Systems},
    VOLUME = {9},
      YEAR = {1989},
    NUMBER = {1},
     PAGES = {101--113},
      ISSN = {0143-3857,1469-4417},
   MRCLASS = {58F11},
  MRNUMBER = {991491},
MRREVIEWER = {Sebastian\ van Strien},
       DOI = {10.1017/S0143385700004831},
       URL = {https://doi.org/10.1017/S0143385700004831},
}

@article {CQ01,
    AUTHOR = {Campbell, James T. and Quas, Anthony N.},
     TITLE = {A generic {$C^1$} expanding map has a singular {S}-{R}-{B}
              measure},
   JOURNAL = {Comm. Math. Phys.},
  FJOURNAL = {Communications in Mathematical Physics},
    VOLUME = {221},
      YEAR = {2001},
    NUMBER = {2},
     PAGES = {335--349},
      ISSN = {0010-3616,1432-0916},
   MRCLASS = {37C40 (28D05 37D20 37D35 37E10)},
  MRNUMBER = {1845327},
MRREVIEWER = {Henk\ Bruin},
       DOI = {10.1007/s002200100491},
       URL = {https://doi.org/10.1007/s002200100491},
}

@article {BT12,
    AUTHOR = {Bruin, Henk and Todd, Mike},
     TITLE = {Transience and thermodynamic formalism for infinitely branched
              interval maps},
   JOURNAL = {J. Lond. Math. Soc. (2)},
  FJOURNAL = {Journal of the London Mathematical Society. Second Series},
    VOLUME = {86},
      YEAR = {2012},
    NUMBER = {1},
     PAGES = {171--194},
      ISSN = {0024-6107,1469-7750},
   MRCLASS = {37D45 (37E05)},
  MRNUMBER = {2959300},
MRREVIEWER = {Pei\ Dong\ Liu},
       DOI = {10.1112/jlms/jdr081},
       URL = {https://doi.org/10.1112/jlms/jdr081},
}

@article {SV97,
    AUTHOR = {Stratmann, B. and Vogt, R.},
     TITLE = {Fractal dimensions for dissipative sets},
   JOURNAL = {Nonlinearity},
  FJOURNAL = {Nonlinearity},
    VOLUME = {10},
      YEAR = {1997},
    NUMBER = {2},
     PAGES = {565--577},
      ISSN = {0951-7715,1361-6544},
   MRCLASS = {28D20 (58F11)},
  MRNUMBER = {1438267},
MRREVIEWER = {Pertti\ Mattila},
       DOI = {10.1088/0951-7715/10/2/014},
       URL = {https://doi.org/10.1088/0951-7715/10/2/014},
}

@article {KN95,
    AUTHOR = {Keller, Gerhard and Nowicki, Tomasz},
     TITLE = {Fibonacci maps re(a{{$\ell$}})visited},
   JOURNAL = {Ergodic Theory Dynam. Systems},
  FJOURNAL = {Ergodic Theory and Dynamical Systems},
    VOLUME = {15},
      YEAR = {1995},
    NUMBER = {1},
     PAGES = {99--120},
      ISSN = {0143-3857,1469-4417},
   MRCLASS = {58F12 (58F03)},
  MRNUMBER = {1314971},
MRREVIEWER = {Remo\ Badii},
       DOI = {10.1017/S0143385700008269},
       URL = {https://doi.org/10.1017/S0143385700008269},
}

@article {Bug85,
    AUTHOR = {Bugiel, Piotr},
     TITLE = {A note on invariant measures for {M}arkov maps of an interval},
   JOURNAL = {Z. Wahrsch. Verw. Gebiete},
  FJOURNAL = {Zeitschrift f\"ur Wahrscheinlichkeitstheorie und Verwandte
              Gebiete},
    VOLUME = {70},
      YEAR = {1985},
    NUMBER = {3},
     PAGES = {345--349},
      ISSN = {0044-3719},
   MRCLASS = {28D05},
  MRNUMBER = {803676},
MRREVIEWER = {Beloslav\ Rie\v can},
       DOI = {10.1007/BF00534867},
       URL = {https://doi.org/10.1007/BF00534867},
}

@article {Ryc76,
    AUTHOR = {Rychlik, Marek},
     TITLE = {Bounded variation and invariant measures},
   JOURNAL = {Studia Math.},
  FJOURNAL = {Polska Akademia Nauk. Instytut Matematyczny. Studia
              Mathematica},
    VOLUME = {76},
      YEAR = {1983},
    NUMBER = {1},
     PAGES = {69--80},
      ISSN = {0039-3223,1730-6337},
   MRCLASS = {28D05 (58F19)},
  MRNUMBER = {728198},
MRREVIEWER = {Nathaniel\ F. G. Martin},
       DOI = {10.4064/sm-76-1-69-80},
       URL = {https://doi.org/10.4064/sm-76-1-69-80},
}

@incollection {MS14,
    AUTHOR = {Moreira, Carlos Gustavo and Smania, Daniel},
     TITLE = {Metric stability for random walks (with applications in
              renormalization theory)},
 BOOKTITLE = {Frontiers in complex dynamics},
    SERIES = {Princeton Math. Ser.},
    VOLUME = {51},
     PAGES = {261--322},
 PUBLISHER = {Princeton Univ. Press, Princeton, NJ},
      YEAR = {2014},
      ISBN = {978-0-691-15929-4},
   MRCLASS = {37A50 (37E05 37E20 37F35 60J10 82B28)},
  MRNUMBER = {3289914},
MRREVIEWER = {Pawe\l\ G\'ora},
}

@InProceedings{Adl73,
author="Adler, Roy L.",
editor="Beck, Anatole",
title="$F$-expansions revisited",
booktitle="Recent Advances in Topological Dynamics",
year="1973",
publisher="Springer Berlin Heidelberg",
address="Berlin, Heidelberg",
pages="1--5",
}

@unpublished{Gun26,
      title={Smooth Circle Covering with a Physical Measure on a Hyperbolic Repelling Fixed Point}, 
      author={Rubio Gunawan},
      year={2026},
      eprint={2602.00293},
      archivePrefix={arXiv},
      primaryClass={math.DS},
      url={https://arxiv.org/abs/2602.00293}, 
      note = {arXiv:2602.00293 [math.DS]}
}
